\documentclass[12pt]{article}

\usepackage{amssymb}
\usepackage{amsmath}

\usepackage[all]{xy}
\input{xy}
\xyoption{all}

\usepackage{geometry}
\newcommand{\N}{\scriptscriptstyle{N}}

\newcommand{\Ll}{\scriptscriptstyle{L}}
\newcommand{\R}{\scriptscriptstyle{R}}

\newcommand{\Mnd}{\mathbf{Mnd}}
\newcommand{\CoMnd}{\mathbf{CoMnd}}

\newcommand{\Hh}{\scriptscriptstyle{H}}

\newcommand{\Ss}{\scriptscriptstyle{S}}
\newcommand{\sS}{\scriptscriptstyle{S}}
\newcommand{\T}{\scriptscriptstyle{T}}
\newcommand{\E}{\scriptscriptstyle{E}}
\newcommand{\K}{\scriptscriptstyle{K}}
\newcommand{\G}{\scriptscriptstyle{G}}
\newcommand{\Pp}{\scriptscriptstyle{P}}
\newcommand{\Q}{\scriptscriptstyle{Q}}

\newcommand{\OpMon}{\mathbf{OpMon}}
\newcommand{\Adj}{\mathbf{Adj}}

\title{Monad and Comonad Objects through \\ 2-adjunctions of the type  Adj-Mnd}

\author{
Adri\'an V\'azquez-M\'arquez\footnote{Contact: avazquez@uiwbajio.mx}\\
Universidad Incarnate Word \\
Campus Baj\'io
}

\date{}

\begin{document}

\maketitle

\begin{center}
\section*{Abstract}
\end{center}

\begin{table}[h]
\centering
\begin{tabular}{p{0.75\textwidth}}
\hline \\
\emph{In this article, the author analyses distributive and mixed
distributive laws and some of their equivalences through the use of 2-adjunctions of
the type $\Adj$-$\Mnd$.}\\
\emph{
  Expressing distributive laws as monad objects in the 2-category $\Mnd(Cat)$ enables one to make an equivalence, through a 2-adjunction, between these distributive laws and monad objects in the 2-category $\Adj_{\R}(Cat)$. These last monad objects thus correspond to monad liftings to algebras of Eilenberg-Moore. Second, the equivalence between these structures and a pair consisting of an Eilenberg-Moore lifting and a Kleisli extension is also analysed. This last result is due to E. Manes and P. Mulry (2010) and it can be restated with an additional naturality on the involved monads.}\\
\emph{
The dual situation, using mixed distributive laws can be analysed too.
}\\
\emph{The objective of this article is based on the use of 2-adjunctions to analyse classic monad theory in order to provide clarity in the proofs and
naturality on the equivalences.}

\\
\hline
\end{tabular}
\end{table}

\noindent 2010 MSC: 18A40, 18C15, 18C20, 18D05, 18D10,
18D35.\\
Keywords: Monads, monoidal categories, adjunctions, 2-categories, 2-adjunctions, Kleisli and Eilenberg-Moore objects.


\section{Introduction}

In \cite{lojj_reke}, the authors use a pair of 2-adjunctions and applied them
to some classical equivalences in monad theory. The reader should check
\cite{cljs_klem} for the construction of such 2-adjunctions.\\ 

The philosophy in \cite{lojj_reke} is that of a revision of equivalences and
bijections in classical monad theory
and analyse them as consequences of higher categorical structures, in this case a 2-adjunction. At first, this approach might
seem more complex than the usual one, but at the expense of an initial complexity,
clarity is gained and not only that but the naturality of the equivalences is obtained at no
cost. Even more, the same context of 2-adjunctions might serve for several applications. \\

The most representative example of this kind of application was the revision
of the equivalence of the lifting of a
monoidal structure to the category of Eilenberg-Moore algebras in terms of an
opmonoidal monad, according to I. Moerdijk \cite{moie_motc}. In that installment, cf.\cite{lojj_reke},
the dual case of the equivalence between an extension of a monoidal
structure over the Kleisli category was analysed too. The reader is invited to check \cite{zama_fotm} in order to follow the approach of higher categorical structures explaining lower structures. \\

In this small article, one of the authors continues with the philosophy in \cite{lojj_reke}, applying the structure of 2-adjunctions of the type $\Adj$-$\Mnd$ to the case of distributive laws and mixed distributive laws. Two seemingly distinct approaches are done for each case. The structure of the article goes as follows. \\

In Section 2.1, the first revision is done and it has to do with the observation due to R. Street in \cite{stro_fotm}, that a distributive law can be seen as a monad object in the 2-category of $\Mnd(Cat)$. Therefore, the 2-adjunction of Eilenberg-Moore serves to transport this monad structure to the 2-category of
$\Adj_{\R}(Cat)$. That is to say, a lifting of the monad structure on the category of Eilenberg-Moore algebras. In this revision, only the naturality on one monad for the equivalence is achived.\\

In Section 2.2, the second revision is done and it has to do with the equivalence of a distributive law with a simultaneous lifting to an Eilenberg-Moore category and an extension to a Kleisli category. This is done by taking, in a simultaneously manner, two 2-adjunctions, one corresponding to the Eilenberg-Moore lifting and the other corresponding to the Kleisli extension. In this revision, the naturality on the two monads for the equivalence is achived.\\

In Section 3.1, the previous plan is applied to the mixed case. The first revision involves the point of view for a mixed distributive law as a comonad object in $\Mnd(Cat)$ and then traslate it to the 2-category $\Adj_{\R}(Cat)$ in order to obtain a lifting of the comonad to the category of Eilenberg-Moore. In this
revision, only the naturality on the monad is found.\\

In Section 3.2, the second approach is done and relies on the equivalence of mixed distributive law with two simultaneous liftings, one to the category of Eilenberg-Moore algebras for the monad and the other to the category of Eilenberg-Moore coalgebras for the comonad. In order to accomplish this, two 2-adjunctions of Eileberg-Moore, of the type $\Adj$-$\Mnd$, have to be handled. In this revision, the naturality on the monad and comonad for this equivalence is found.\\

The notation and conventions for this article goes as follows. We take the direction an adjunction as the direction of its corresponding left adjoint functor. We allow ourselves to take the notation of a monad $(S,\mu^{\Ss},\eta^{\Ss})$ on the category $\mathcal{C}$ as $(\mathcal{C}, S, \mu^{\Ss},\eta^{\Ss})$ or $(\mathcal{C},S)$, in a compact form, where the multiplication and unit are left undestood. While speaking of $n$-cells, and most of the time, only the arrow corresponding to the $n$-$1$-cells will be written and the $r$-cells below will be left understood.\\ 

As far as notation is concerned, for the free-forgetful adjunction
of the Eilenberg-Moore category for the monad $(S, \mu^{\Ss},\eta^{\Ss})$, the notation
$F^{\Ss}\dashv U^{\Ss}: \mathcal{C}\longrightarrow \mathcal{C}^{\Ss}$ is used. In
turn, the Kleisli adjunction is denotated as $D_{\Ss}\dashv V_{\Ss}:\mathcal{C}\longrightarrow\mathcal{C}_{\Ss}$.
Given an adjunction $L\dashv R: \mathcal{C}\longrightarrow \mathcal{D}$, the
unit will be denoted as $\eta^{\R\Ll}$ and the counit as
$\varepsilon^{\Ll\R}$. This notation is used for the sake of avoiding multiple different greek letters to denote units and counits. \\

In the whole of the present article $Cat$ will stand for the 2-category of small categories, functors and natural transformations.\\

An appendex is provided in order to explain, up to a certain detail, the construction of
the 2-adjunctions used in this article.\\

\section{Distributive Laws and Monad Objects}

We begin this section by giving the definition of a distributive law.

\newtheorem{1705141458}{Definition}[section]
\begin{1705141458}

A distributive law from the monad $(\mathcal{C}, S)$ to the monad
$(\mathcal{C}, T)$ is a natural transformation $\varphi:ST\longrightarrow TS$
such that the following diagrams commute:

\begin{equation}\label{1510071248}
\begin{array}{ccc}
\xy<1cm,0cm>;
(0,24) *+{SST} = "a11",
(30,24) *+{STS} = "a12",
(60,24) *+{TSS} = "a13",
(0,0) *+{ST} = "a21",
(60,0) *+{TS} = "a23",
\POS "a11" \ar^{S\varphi} "a12",
\POS "a12" \ar^{\varphi S} "a13",
\POS "a13" \ar^{T\mu^{\Ss}} "a23",
\POS "a11" \ar_{\mu^{\Ss}T} "a21",
\POS "a21" \ar_{\varphi} "a23",
\endxy && 
\xy<1cm,0cm>;
(16,16) *+{T} = "a12",
(0,0) *+{ST} = "a21",
(32,0) *+{TS} = "a23",
\POS "a12" \ar^{T\eta^{\Ss}} "a23",
\POS "a12" \ar_{\eta^{\Ss}T} "a21",
\POS "a21" \ar_{\varphi} "a23",
\endxy
\end{array}
\end{equation}

\begin{equation}\label{1510071249}
\begin{array}{ccc}
\xy<1cm,0cm>;
(0,24) *+{STT} = "a11",
(30,24) *+{TST} = "a12",
(60,24) *+{TTS} = "a13",
(0,0) *+{ST} = "a21",
(60,0) *+{TS} = "a23",
\POS "a11" \ar^{\varphi T} "a12",
\POS "a12" \ar^{T\varphi} "a13",
\POS "a13" \ar^{\mu^{\T}S} "a23",
\POS "a11" \ar_{S\mu^{\T}} "a21",
\POS "a21" \ar_{\varphi} "a23",
\endxy && 
\xy<1cm,0cm>;
(16,16) *+{S} = "a12",
(0,0) *+{ST} = "a21",
(32,0) *+{TS} = "a23",
\POS "a12" \ar^{\eta^{\T}S} "a23",
\POS "a12" \ar_{S\eta^{\T}} "a21",
\POS "a21" \ar_{\varphi} "a23",
\endxy
\end{array}
\end{equation}

\end{1705141458}

\subsection{Monad Objects}

In 1969, J. Beck proved, in \cite{eckb_setc-bejo_dila}, that a distributive law is equivalent to a
lifting of the monad structure $T$ on $\mathcal{C}$ to a monad structure $\widehat{T}$ on the
category of Eilenberg-Moore algebras, namely $\mathcal{C}^{\Ss}$.\\

R. Street, in \cite{stro_fotm}, noticed that a distributive law $\varphi$, from $(S, \mu^{\sS}, \eta^{\sS})$ to $(T, \mu^{\T}, \eta^{\sS})$, can be seen as a monad object in $\Mnd(Cat)$, with endo 1-cell $(T, \varphi)$ over the object $(\mathcal{C}, S)$.\\

The J. Beck's equivalence is taken within the context of a 2-adjunction of the type $\Adj$
- $\Mnd$, cf. \cite{cljs_klem} and \cite{lojj_reke}.\\

\newtheorem{1705201918}{Theorem}[subsection]
\begin{1705201918}\label{1705201918}
There is a bijection between the following structures

\begin{enumerate}

\item Lifting monads $(\widehat{T}, \mu^{\widehat{\T}}, \eta^{\widehat{\T}})$, on $\mathcal{C}^{\sS}$, for $(T, \mu^{\T}, \eta^{\T})$ on $\mathcal{C}$. That is to say, monad objects in $\Adj_{\R}(Cat)$.  

\item Distributive laws from the monad $(S, \mu^{\sS}, \eta^{\sS})$, on $\mathcal{C}$, to the monad $(T, \mu^{\T}, \eta^{\T})$. That is to say, monad objects in $\Mnd(Cat)$ of the form $((\mathcal{C}, S), (T, \varphi), \mu^{\T}, \eta^{\T})$.

\end{enumerate}

This bijection is natural on the monad $(S, \mu^{\sS}, \eta^{\sS})$.

\end{1705201918}

$\phantom{br}$

\noindent \emph{Proof}:\\

Take the 2-adjunction of Eilenberg-Moore $\Phi_{\E} \dashv \Psi_{\E}$

\begin{equation}\label{1510081248}
\begin{array}{cc}
\xymatrix@C=1.6cm{
\Adj_{\R}(Cat) \ar@<-3pt>[r]_{\Phi_{\E}} & \Mnd(Cat)\ar@<-3pt>[l]_{\Psi_{\!\E}} 
} &\!\!\! ,
\end{array}
\end{equation}

$\phantom{br}$

From this 2-adjunction, we can get an isomorphism of categories for any pair of objects $L\dashv R$ in $\Adj_{\R}(Cat)$ and $(\mathcal{Z}, T)$ in $\Mnd(Cat)$

\begin{equation*}
Hom_{\Adj_{\R}}(L\dashv R, \Psi_{\E}(\mathcal{Z}, T)) \cong Hom_{\Mnd}(\Phi_{\E}(L\dashv R), (\mathcal{Z}, T))
\end{equation*}

$\phantom{br}$

If in the previous isomorfism we make  $L\dashv R = D^{\Ss}\dashv U^{\Ss}$ y
$(\mathcal{Z}, T) = (\mathcal{C}, S)$, then 

\begin{equation}\label{1510071407}
Hom_{\Adj_{\R}}(D^{\Ss}\dashv U^{\Ss}, D^{\Ss}\dashv U^{\Ss} ) \cong Hom_{\Mnd}(  (\mathcal{C}, S), (\mathcal{C}, S))
\end{equation}

$\phantom{br}$

This isomorphism takes, for example, 1-cells in $\Adj_{\R}(Cat)$ of the form

\begin{equation*}
\xy<0cm,1cm>;
(0,20) *+{\mathcal{C}} = "a11",
(20,20) *+{\mathcal{C}} = "a12",
(0,0) *+{\phantom{.}\mathcal{C}^{\Ss}} = "a21",
(20,0) *+{\phantom{.}\mathcal{C}^{\Ss}} = "a22",
\POS "a21" \ar^{U^{\Ss}} "a11",
\POS "a11" \ar^{T} "a12",
\POS "a21" \ar_{\widehat{T}} "a22",
\POS "a22" \ar_{U^{\Ss}} "a12",
\endxy
\end{equation*}

\noindent to 1-cells in $\Mnd(Cat)$ of the form $(T, \varphi):(\mathcal{C}, S)\longrightarrow (\mathcal{C}, S)$ such that $\varphi:ST\longrightarrow TS$ is a natural transformation making the following diagrams commutative

\begin{equation*}
\xymatrix@C=1.9cm@R=0.9cm{
SST \ar[r]^-{S\varphi}\ar[dd]_{\mu^{\Ss}\! T} & STS\ar[r]^-{\varphi S} & TSS\ar[dd]^{T\mu^{\Ss}}\\
&T\ar[ld]_{\eta^{\Ss}T}\ar[rd]^{T\eta^{\Ss}}&\\
ST\ar[rr]_{\varphi} & &     TS
} 
\end{equation*}

$\phantom{br}$

Note that the two diagrams were \emph{merged} into one.\\

Consider now the following 2-cells in $\Adj_{\R}(Cat)$:

\begin{equation*}
\begin{array}{cccc}
\xy<1cm,0cm>
\POS (0, 0) *+{\mathcal{C}} = "c1",
\POS (28, 0) *+{\mathcal{C}} = "c2", 
\POS (0, -30) *+{\mathcal{C}^{\Ss}} = "d1",
\POS (28, -30) *+{\mathcal{C}^{\Ss}} = "d2", 
\POS (14, 0) *+{\mu} = "e1",
\POS (14, -30) *+{\widehat{\mu}} = "e2",
\POS (8, 3) = "a",
\POS (8, -2.5) = "b",
\POS (8,-27) = "ap",
\POS (8, -32.5) = "bp",
\POS "c1" \ar@/^1.2pc/^{TT} "c2",
\POS "c1" \ar@/_1.2pc/_{T} "c2",
\POS "c1" \ar@<-2pt>_{D^{\Ss}} "d1",
\POS "d1" \ar@<-2pt>_{U^{\Ss}} "c1",
\POS "c2" \ar@<-2pt>_{D^{\Ss}} "d2",
\POS "d2" \ar@<-2pt>_{U^{\Ss}} "c2",
\POS "d1" \ar@/^1.2pc/^{\widehat{T}\widehat{T}} "d2",
\POS "d1" \ar@/_1.2pc/_{\widehat{T}} "d2",
\POS "a" \ar "b",
\POS "ap" \ar "bp",
\endxy & 
\xy<1cm,0cm>
\POS (0,-31.5) *+{,} 
\endxy
&
\xy<1cm,0cm>
\POS (0, 0) *+{\mathcal{C}} = "c1",
\POS (28, 0) *+{\mathcal{C}} = "c2", 
\POS (0, -30) *+{\mathcal{C}^{\Ss}} = "d1",
\POS (28, -30) *+{\mathcal{C}^{\Ss}} = "d2", 
\POS (14, 0) *+{\eta} = "e1",
\POS (14, -30) *+{\widehat{\eta}} = "e2",
\POS (8, 3) = "a",
\POS (8, -2.5) = "b",
\POS (8,-27) = "ap",
\POS (8, -32.5) = "bp",
\POS "c1" \ar@/^1.2pc/^{1_{\mathcal{C}}} "c2",
\POS "c1" \ar@/_1.2pc/_{T} "c2",
\POS "c1" \ar@<-2pt>_{D^{\Ss}} "d1",
\POS "d1" \ar@<-2pt>_{U^{\Ss}} "c1",
\POS "c2" \ar@<-2pt>_{D^{\Ss}} "d2",
\POS "d2" \ar@<-2pt>_{U^{\Ss}} "c2",
\POS "d1" \ar@/^1.2pc/^{1_{\mathcal{C}^{\Ss}}} "d2",
\POS "d1" \ar@/_1.2pc/_{\widehat{T}} "d2",
\POS "a" \ar "b",
\POS "ap" \ar "bp",
\endxy & 
\xy<1cm,0cm>
\POS (0,-31.5) *+{.} 
\endxy
\end{array} 
\end{equation*}

Therefore, the following conditions are fulfilled

\begin{eqnarray*}
U^{\Ss}\widehat{\mu}^{\T}  &=& \mu^{\T} U^{\Ss}\\
U^{\Ss}\widehat{\eta}^{\T} &=& \eta^{\T} U^{\Ss}
\end{eqnarray*}

Note that the counit of the 2-adjunction $\varepsilon^{\Phi\Psi}$, on any 0-cell $F^{\Ss}\dashv U^{\Ss}$, is the identity $(1_{\mathcal{C}}, 1_{\mathcal{C}^{\Ss}})$. Therefore, the isomorphic image of the 2-cells $(\mu^{\T}, \mu^{\widehat{\T}})$ and $(\eta^{\T}, \eta^{\widehat{\T}})$ is the pair of 2-cells in $\Mnd(Cat)$

\begin{eqnarray*}
\mu^{\T}&:&(T,\varphi)\cdot(T,\varphi)\longrightarrow (T, \varphi):
(\mathcal{C}, S)\longrightarrow (\mathcal{C}, S),\\
\eta^{\T}&:&(1_{\mathcal{C}}, 1_{\Ss})\longrightarrow (T, \varphi):
(\mathcal{C}, S)\longrightarrow (\mathcal{C}, S).
\end{eqnarray*}

$\phantom{br}$

That is to say, they fulfill the commutative of the following diagrams

\begin{equation*}
\xy<1cm,0cm>;
(0,24) *+{STT} = "a11",
(30,24) *+{TST} = "a12",
(60,24) *+{TTS} = "a13",
(0,0) *+{ST} = "a21",
(60,0) *+{TS} = "a23",
\POS "a11" \ar^{\varphi T} "a12",
\POS "a12" \ar^{T\varphi} "a13",
\POS "a13" \ar^{\mu^{\T}S} "a23",
\POS "a11" \ar_{S\mu^{\T}} "a21",
\POS "a21" \ar_{\varphi} "a23",
\endxy
\end{equation*}

and

\begin{equation*}
\xy<1cm,0cm>;
(15,18) *+{S} = "a12",
(0,0) *+{ST} = "a21",
(30,0) *+{TS} = "a23",
\POS "a12" \ar^{\eta^{\T}\!S} "a23",
\POS "a12" \ar_{S\eta^{\T}} "a21",
\POS "a21" \ar_{\varphi} "a23",
\endxy
\end{equation*}

$\phantom{br}$

In summary, the isomorphic image for a monad $\big( (T, \widehat{T}), (\mu^{\T}, \widehat{\mu}^{\T}), (\eta^{\T}, \widehat{\eta}^{\T})\big)$, on $F^{\Ss}\dashv U^{\Ss}$, is $\big((\mathcal{C},S), (T,\varphi ), \mu^{\T},\eta^{\T}\big)$. Where, the natural transformation $\varphi$ has the following form

\begin{equation}\label{1510081202}
\varphi = U^{\Ss}\lambda_{\T} = U^{\Ss}\varepsilon^{FU}\widehat{T}F^{\Ss}\circ U^{\Ss}F^{\Ss}T\eta^{UF}.
\end{equation}

$\phantom{br}$

\noindent according to the construction of the 2-functor $\Phi_{\E}$, see Appendex.\\

Due to the isomorphism of categories given by \eqref{1510071407}, the previous bijection is natural on the monad $(\mathcal{C}, S, \mu^{\Ss}, \eta^{\Ss})$.

\begin{flushright}
$\Box$
\end{flushright}

Note that any 2-functor sends any monad object to a monad object, but the equivalence found here is a particular one.\\

\subsection{Combining Eilenberg-Moore liftings and Kleisli Extensions}

The naturality on both monads can be finally stablished but
with the addition of a second 2-adjunction. In order to do so, consider the following theorem due to E. Manes \& P. Mulry, \cite{maep_mocg}.

\newtheorem{1501131518}{Theorem}[subsection]
\begin{1501131518}\label{1501131518}

Given the following monads $(S, \mu^{\Ss}, \eta^{\Ss})$ \& $(T, \mu^{\T},
\eta^{\T})$ on $\mathcal{C}$. A natural transformation $\varphi: ST\longrightarrow TS$ is a
 distributive law from $S$ to $T$ iff classifies both as a lifting of
Eilenberg-Moore, for the endofunctor $T$, and as a Kleisli Extension, for the
endofunctor $S$, according to the following diagrams

\begin{equation}\label{1510081116}
\begin{array}{ccc}
\xy<0cm,1cm>;
(0,22) *+{\mathcal{C}} = "a11",
(22,22) *+{\mathcal{C}} = "a12",
(0,0) *+{\phantom{.}\mathcal{C}^{\Ss}} = "a21",
(22,0) *+{\phantom{.}\mathcal{C}^{\Ss}} = "a22",
\POS "a21" \ar^{U^{\Ss}} "a11"
\POS "a11" \ar^{T} "a12",
\POS "a21" \ar_{\widehat{T}} "a22",
\POS "a22" \ar_{U^{\Ss}} "a12",
\endxy
& &
\xy<0cm,1cm>;
(0,22) *+{\mathcal{C}} = "a11",
(22,22) *+{\mathcal{C}} = "a12",
(0,0) *+{\phantom{.}\mathcal{C}_{\T}} = "a21",
(22,0) *+{\phantom{.}\mathcal{C}_{\T}} = "a22",
\POS "a11" \ar_{D_{\T}} "a21"
\POS "a11" \ar^{S} "a12",
\POS "a21" \ar_{\widetilde{S}} "a22",
\POS "a12" \ar^{D_{\T}} "a22",
\endxy
\end{array}
\end{equation}

\end{1501131518}

\emph{Proof}:\\

Half of the proof has been given, that is to say, there is a bijective correspondence between liftings of the form $(T, \widehat{T}):F^{\sS}\dashv U^{\sS}\longrightarrow F^{\sS}\dashv U^{\sS}$ and morphisms of monads of the form $(T, \varphi): (\mathcal{C}, S)\longrightarrow(\mathcal{C}, S)$.\\

In order to establish a similar bijective correspondence, natural in the monad $(\mathcal{C}, T, \mu^{\T}, \eta^{\T})$, another 2-adjunction of the type $\Adj$-$\Mnd$ has to be considered, cf. \cite{cljs_klem} and \cite{lojj_reke}. The so-called Kleisli 2-adjunction:

\begin{equation*}
\begin{array}{cc}
\xymatrix@C=1.6cm{
\Mnd^{\bullet}(Cat)  \ar@<-3pt>[r]_{\Psi_{\K}} & \Adj_{\Ll}(Cat)\ar@<-3pt>[l]_{\Phi_{\!\K}} 
} &\!\!\! ,
\end{array}
\end{equation*}

$\phantom{br}$

The construction of such a 2-adjunction is described in the Appendex.\\

Due to the 2-adjunction, there exist the following isomorphism of categories, natural on the monad $(\mathcal{C}, T)$, in $\Mnd^{\bullet}(Cat)$, and on the adjunction $D_{\T}\dashv V_{\T}$, in $\Adj_{\Ll}(Cat)$,

\begin{equation*}
Hom_{\Mnd^{\bullet}}((\mathcal{C}, T), \Phi_{\K}(D_{\T}\dashv V_{\T}))
\cong Hom_{\Adj_{\Ll}}(\Psi_{\K}(\mathcal{C}, T), D_{\T}\dashv V_{\T})
\end{equation*}

$\phantom{br}$

It is an isomorphism of the form 

\begin{equation*}
Hom_{\Mnd^{\bullet}}((\mathcal{C}, T), (\mathcal{C}, T))
\cong Hom_{\Adj_{\Ll}}(D_{\T}\dashv V_{\T}, D_{\T}\dashv V_{\T})
\end{equation*}

$\phantom{br}$

Therefore, extensions to the Kleisli Category  $\mathcal{C}_{\T}$, that is to say, objects in the category $Hom_{\Adj_{\Ll}}(D_{\T}\dashv V_{\T}, D_{\T}\dashv V_{\T})$,

\begin{equation*}
\xy<0cm,1cm>;
(0,22) *+{\mathcal{C}} = "a11",
(22,22) *+{\mathcal{C}} = "a12",
(0,0) *+{\phantom{.}\mathcal{C}_{\T}} = "a21",
(22,0) *+{\phantom{.}\mathcal{C}_{\T}} = "a22",
\POS "a11" \ar_{D_{\T}} "a21",
\POS "a11" \ar^{S} "a12",
\POS "a21" \ar_{\widetilde{S}} "a22",
\POS "a12" \ar^{D_{\T}} "a22",
\endxy
\end{equation*}

\noindent are in bijection with natural transformations  $\varphi:ST\longrightarrow TS$ such that the following diagrams commute 

\begin{equation*}
\xymatrix@C=1.9cm@R=0.9cm{
STT \ar[r]^-{\varphi T}\ar[dd]_{S\mu^{\T}} & TST\ar[r]^-{T\varphi} & TTS\ar[dd]^{\mu^{\T}S}\\
&S\ar[ld]_{S\eta^{\T}}\ar[rd]^{\eta^{\T}S}&\\
ST\ar[rr]_{\varphi} & &     TS
} 
\end{equation*}

$\phantom{br}$

These commutative diagrams are precisely the ones in \eqref{1510071249} completing a distributive law from $S$ to $T$.\\

The only thing that has to be taken into account is the following. If one starts from a lifting over the Eilenberg-Moore category and an extension over the
Kleisli category, the induced natural transformation $\varphi: ST\longrightarrow TS$ have to be the same. This requirement can be achived by asking the following equality to hold

\begin{equation*}
U^{\Ss}\lambda_{\T} = \rho_{\Ss}D_{\T},
\end{equation*}

$\phantom{br}$

\noindent where $\lambda_{\T}$ is the mate of the first commutative diagram in \eqref{1510081116} and $\rho_{\T}$ is the mate of the second commutative diagram.\\

The previous requirement can be expressed also by 

\begin{equation*}
\Phi_{\E}(T,\widehat{T}, \lambda_{\T}) = \Psi_{\K}(S, \widetilde{S}, \rho_{\Ss})
\end{equation*}

\begin{flushright}
$\Box$
\end{flushright}

$\phantom{br}$

At this point Theorem \ref{1501131518} can be restated with an additional naturality on the monads $(S,\mu^{\Ss},\eta^{\Ss})$ and  $(T, \mu^{\T},\eta^{\T})$ on $\mathcal{C}$.

\section{Mixed Distributive Laws}

The previous procedure can be applied to the case of mixed distributive laws.

\newtheorem{1412081737}{Definition}[section]
\begin{1412081737}

A mixed distributive law from the monad $(S, \mu^{\Ss}, \eta^{\Ss})$ to the
comonad $(G, \delta^{\G}, \varepsilon^{\G})$, on $\mathcal{C}$, is a natural transformation
$\psi:SG\longrightarrow GS$ such that the following diagrams commute

\begin{equation}\label{1510081240}
\begin{array}{ccc}
\xy<1cm,0cm>;
(0,24) *+{SSG} = "a11",
(30,24) *+{SGS} = "a12",
(60,24) *+{GSS} = "a13",
(0,0) *+{SG} = "a21",
(60,0) *+{GS} = "a23",
\POS "a11" \ar^{S\psi} "a12",
\POS "a12" \ar^{\psi S} "a13",
\POS "a13" \ar^{G\mu^{\Ss}} "a23",
\POS "a11" \ar_{\mu^{\Ss}G} "a21",
\POS "a21" \ar_{\psi} "a23",
\endxy && 
\xy<1cm,0cm>;
(15,18) *+{G} = "a12",
(0,0) *+{SG} = "a21",
(30,0) *+{GS} = "a23",
\POS "a12" \ar^{G\eta^{\Ss}} "a23",
\POS "a12" \ar_{\eta^{\Ss}\! G} "a21",
\POS "a21" \ar_{\psi} "a23",
\endxy
\end{array}
\end{equation}

\begin{equation}\label{1510081242}
\begin{array}{ccc}
\xy<1cm,0cm>;
(0,24) *+{SGG} = "a11",
(30,24) *+{GSG} = "a12",
(60,24) *+{GGS} = "a13",
(0,0) *+{SG} = "a21",
(60,0) *+{GS} = "a23",
\POS "a11" \ar^{\psi G} "a12",
\POS "a12" \ar^{G\psi } "a13",
\POS "a23" \ar_{\delta^{\G}S} "a13",
\POS "a21" \ar^{S\delta^{\G}} "a11",
\POS "a21" \ar_{\psi} "a23",
\endxy && 
\xy<1cm,0cm>;
(15,18) *+{S} = "a12",
(0,0) *+{SG} = "a21",
(30,0) *+{GS} = "a23",
\POS "a23" \ar_{\varepsilon^{\G}\! S} "a12",
\POS "a21" \ar^{S\varepsilon^{\G}} "a12",
\POS "a21" \ar_{\psi} "a23",
\endxy
\end{array}
\end{equation}

\end{1412081737}

\subsection{Comonad Objects}

In the same way as we did for the case of monads objects in Section 2.1.

\newtheorem{1705222141}{Theorem}[subsection]
\begin{1705222141}
There is bijection between the following structures:

\begin{enumerate}

\item Mixed distributive laws from the monad $(S, \mu^{\sS}, \eta^{\sS})$, on $\mathcal{C}$, to the comonad $(G, \delta^{\G}, \varepsilon^{\G})$. That is to say, comonad objects in $\Mnd(Cat)$.

\item Lifted comonads $(\widehat{G}, \delta^{\widehat{\G}}, \varepsilon^{\widehat{\G}})$, on $\mathcal{C}^{\sS}$, for $(G, \delta^{\G}, \varepsilon^{\G})$ on $\mathcal{C}$. That is to say, comonad objects in $\Adj_{\R}(Cat)$.\\

This bijection is natural on the monad $(S, \mu^{\sS}, \eta^{\sS})$, on $\mathcal{C}$.
  
\end{enumerate}

\end{1705222141}

\begin{flushright}
$\Box$
\end{flushright}

\subsection{Combining Eilenberg-Moore liftings of algebras and coalgebras}

$\phantom{br}$

A similar theorem for mixed distributive laws can be stated that involve equivalences of certain structures. In the following case through a pair of Eilenberg-Moore liftings.\\

\newtheorem{1501140312}{Theorem}[subsection]
\begin{1501140312}\label{1501140312}

Given the monad $(S, \mu^{\Ss}, \eta^{\Ss})$ and the comonad $(G, \delta^{\G},
\varepsilon^{\G})$ on $\mathcal{C}$. A natural transformation $\psi:
SG\longrightarrow GS$ is a mixed distributive law from $S$ to $G$ iff the
endofunctors classify both as Eilenberg-Moore liftings:

\begin{equation}\label{1510081654}
\begin{array}{ccc}
\xy<0cm,1cm>;
(0,22) *+{\mathcal{C}} = "a11",
(22,22) *+{\mathcal{C}} = "a12",
(0,0) *+{\phantom{.}\mathcal{C}^{\Ss}} = "a21",
(22,0) *+{\phantom{.}\mathcal{C}^{\Ss}} = "a22",
\POS "a21" \ar^{U^{\Ss}} "a11"
\POS "a11" \ar^{G} "a12",
\POS "a21" \ar_{\widehat{G}} "a22",
\POS "a22" \ar_{U^{\Ss}} "a12",
\endxy
& &
\xy<0cm,1cm>;
(0,22) *+{\mathcal{C}^{\G}} = "a11",
(22,22) *+{\mathcal{C}^{\G}} = "a12",
(0,0) *+{\mathcal{C}} = "a21",
(22,0) *+{\mathcal{C}} = "a22",
\POS "a11" \ar_{U^{\G}} "a21"
\POS "a11" \ar^{\widehat{S}} "a12",
\POS "a21" \ar_{S} "a22",
\POS "a12" \ar^{U^{\G}} "a22",
\endxy
\end{array}
\end{equation}

\end{1501140312}

$\phantom{br}$

Note the change of direction for the Eilenberg-Moore coalgebra lifting. This
is due to the convention of taking the direction of the left adjoint functor.\\

\noindent \emph{Proof}:\\

The structure of the proof can be stablished as follows. First, we relate the diagrams in \eqref{1510081240} with the lifting \eqref{1510081654}. Second, we relate the diagrams in \eqref{1510081242} with the second lifting diagram in \eqref{1510081654}.\\

The first part was already done in the proof of \ref{1705201918} just change the endofunctor $T$ by $G$.\\

For the second part, we have to consider another 2-adjunction, that of Eilenberg-Moore for comonads. To begin with, we take the \emph{co-opossite} 2-category $Cat^{co}$. We note the following isomorphisms

\begin{equation*}
\Adj_{\R}(Cat^{co})\cong \Adj_{\Ll}(Cat)\:,\:\: \Mnd(Cat^{co})\cong \CoMnd(Cat)
\end{equation*}

$\phantom{br}$

If the 2-category $\CoMnd(Cat)$ admits the construction of algebras then the following 2-adjunction can be constructed

\begin{equation*}
\begin{array}{cc}
\xymatrix@C=1.6cm{
\Adj_{\Ll}(Cat) \ar@<-3pt>[r]_{\vec{\Phi}_{\E}} & \CoMnd(Cat)\ar@<-3pt>[l]_{\vec{\Psi}_{\!\E}} 
} &\!\!\! ,
\end{array}
\end{equation*}

$\phantom{br}$

From this 2-adjunction, and for the 0-cells $L\dashv R$ in $Adj_{\Ll}(Cat)$ and $(\mathcal{Z}, G)$ comonad in $\CoMnd(Cat)$, there exist a natural isomorphism 

\begin{equation*}
Hom_{\Adj_{\Ll}}(L\dashv R, \vec{\Psi}_{\E}(\mathcal{Z}, G)) \cong Hom_{\CoMnd}(\vec{\Phi}_{\E}(L\dashv R), (\mathcal{Z}, G))
\end{equation*}

$\phantom{br}$

In particular, if $L\dashv R = U^{\G}\dashv D^{\G}$, the Eileberg-Moore adjunction for the corresponding comonad, and $\mathcal{Z} = \mathcal{C}$

\begin{equation*}
Hom_{\Adj_{\Ll}}(U^{\G}\dashv D^{\G}, U^{\G}\dashv D^{\G}) \cong Hom_{\CoMnd}((\mathcal{C}, G), (\mathcal{C}, G))
\end{equation*}

$\phantom{br}$

Therefore, there is a bijective correspondence between objects in the category $Hom_{\Adj_{\Ll}}(U^{\G}\dashv D^{\G}, U^{\G}\dashv D^{\G})$, that is to say, liftings to the Eilenberg-Moore category of coalgebras as in \eqref{1510081654} and objects in the category  $Hom_{\CoMnd}((\mathcal{C}, G), (\mathcal{C}, G))$,
that is to say, pairs of commutative diagrams as in \eqref{1510081242}. \\

Finally, we have to consider the compatibility of the following image for the pair of 2-adjunctions

\begin{equation*}
\Phi_{\E}(G, \widehat{G}, \lambda_{\G}) = \vec{\Phi}_{\E}(S, \widehat{S}, \rho_{\Ss})
\end{equation*}

\noindent or 

\begin{equation*}
U^{\Ss}\lambda_{\G} = U^{\G}\!\rho_{\Ss}
\end{equation*}

Note that $\rho_{\sS}$ is the mate of the second commutative natural transformation in \eqref{1510081654}.\\

The previous equivalence is natural on the monad $(\mathcal{C}, S, \mu^{\Ss}, \eta^{\Ss})$
and on the comonad $(\mathcal{C}, G, \delta^{\G}, \varepsilon^{\G})$.

\begin{flushright}
$\Box$
\end{flushright}

\section{Formal set up}

All the discussion can be taken into the the formal set up without any trouble, that is to say, a general 2-category $\mathcal{A}$ can be considered with the requirement that it, and its duals $\mathcal{A}^{op}$ and $\mathcal{A}^{co}$, accepts the construction of algebras. The reader is compelled to check \cite{lojj_reke} for the construction of the Kleisli 2-adjunction in this general set up ($\mathcal{A}$ and $\mathcal{A}^{op}$) and to check \cite{stro_fotm} for what involves for the 2-category $\mathcal{A}$, and its duals, to accept the construction of algebras. This last reference provides also the necessary subjacent calculations for the construction of the Eilenberg-Moore 2-adjunction for a general 2-category $\mathcal{A}$.

\section{Epilogue}

The 2-adjunction of Eileberg-Moore can be used also in the context of monoidal categories as follows. Let $\big((S, \psi^{\Ss}), \mu^{\Ss}, \eta^{\Ss}\big)$ be a opmonoidal monad on the monoidal category $(\mathcal{C},\otimes, I)$. In 2002,  I. Moerdijk induced a monoidal structure on the category of Eilenberg-Moore algebras, $\mathcal{C}^{\Ss}$, for such an opmonoidal monad, cf. \cite{moie_motc}. In 2010, M. Zawadoski proved that the 2-category of opmonoidal monads accepts the construction of algebras, cf. \cite{zama_fotm}. The Eilenberg-Moore object for an opmonoidal monad is precisely the induced monoidal structure of I. Moerdijk, therefore the following 2-adjunction can be constructed

\begin{equation}\label{1511052041}
\begin{array}{cc}
\xymatrix@C=1.6cm{
\Adj_{\R}(\OpMon(Cat)) \ar@<-3pt>[r]_{\Phi_{\E}} & \Mnd(\OpMon(Cat))\ar@<-3pt>[l]_{\Psi_{\!\E}} 
} &\!\!\! ,
\end{array}
\end{equation}

$\phantom{br}$

Note that the notation for this 2-category, $\OpMon(Cat)$, is different from the one in \cite{zama_fotm}, namely $\mathbf{Mon}_{op}(Cat)$.\\

Therefore for a 0-cell $Adj_{\R}(\OpMon(Cat))$, \emph{i.e.} an opmonoidal adjunction $(L, \tau^{\Ll})\dashv (R,\tau^{\R})$, cf. \cite{bral_homm}, for the precise details; and for a 0-cell in $\Mnd(\OpMon(Cat))$, $\big((\mathcal{C}, \otimes), (S, \psi^{\Ss}),\mu^{\Ss}, \eta^{\Ss} \big)$, an opmonoidal monad, there is an isomorphism of categories as follows:

\setlength\arraycolsep{0.6pt}
\begin{eqnarray*}
Hom_{\Adj_{\R}(\OpMon(Cat))}\Big((L,\tau^{\Ll})\dashv (R,\tau^{\R}) &,& \Psi_{\E}\big((\mathcal{C}, \otimes), (S,
\psi^{\Ss}),\mu^{\Ss}, \eta^{\Ss}\big)\Big)\\
&\cong& \\
Hom_{\Mnd(\OpMon(Cat))}\Big(\Phi_{\E}\big((L,
\tau^{\Ll})\dashv (R,\tau^{\R})\big) &,& ((\mathcal{C}, \otimes), (S,
\psi^{\Ss}),\mu^{\Ss}, \eta^{\Ss}))\Big)
\end{eqnarray*}

This isomorphism may help to set up a 2-categorical approach in the characterization of Hopf monads using opmonoidal adjunctions, cf. \cite{agmc_gehm} and \cite{bral_homm}. 

\section{Discussion}

In \cite{eckb_setc-bejo_dila}, the author Jon Beck did find another equivalence for the characterisation of distributive laws in terms of the product monad. The author of this article did not get this other equivalence in terms of the 2-adjunction. This might use other properties of the 2-adjunction.\\

The author consider that the present approach is very good to check classical results in monad
theory due to the fact that adds clarity to the proofs. On the account of some
complexity, at the beginning, the explanations are smoother and the
applications can be several. Note that the same 2-adjunction can be reused to
produce different equivalences against a usual proof where the procedure
might not be reused in an straightforward way.\\

\section{Conclusions and Future Work}

These examples can be added to the list started in \cite{lojj_reke}, which
accounts for the applicability of this procedure.\\

As fas as the future work is concerned, it remains to explore the 2-adjunction in order to get the equivalence for the product monad.\\

Finally, the author encourages the readers to get other examples that fit into
this context in order to understand what the higher category theory remains to
say about classical monad theory or to get examples in category
theory, logic or functional programming where the 2-adjunction context can be applied.

\section{Appendex}

In this appendex, we list the three 2-adjunctions that took a role in this
article, with a two-fold
purpose. First, the purpose to have them at hand in case someone wants to
apply them. Second, to explain its use within this article but without
interfering to the continuity of the reading. 

\subsection*{The Eilenberg-Moore 2-adjunction for monads}

The Eilenberg-Moore adjunction for monads 

\begin{equation}
\begin{array}{cc}
\xymatrix@C=1.6cm{
\Adj_{\R}(Cat) \ar@<-3pt>[r]_{\Phi_{\E}} & \Mnd(Cat)\ar@<-3pt>[l]_{\Psi_{\!\E}} 
} &\!\!\! ,
\end{array}
\end{equation}

\noindent is comprised of several structures. To begin with, the 2-category $\Adj_{\R}(Cat)$ with the following $n$-cells

\begin{enumerate}

\item 0-cells are adjunctions $L\dashv R: \mathcal{C}\longrightarrow\mathcal{X}$ and $\overline{L}\dashv\overline{R}:\mathcal{D}\longrightarrow\mathcal{Y}$.

\item 1-cells are pairs of functors and a natural transformation $(J,K,\lambda):L\dashv R\longrightarrow \overline{L}\dashv\overline{R}$. Where $J:\mathcal{C}\longrightarrow\mathcal{D}$ and $K:\mathcal{X}\longrightarrow\mathcal{Y}$ and the square involving the right adjoint funtors commute, $KR = \overline{R}J$, hence the name after this commutativity. The natural transformation $\lambda$ corresponds to the mate of the commutative diagram.

\item 2-cells are pairs of natural transformations $(\alpha,
  \beta):(J,K,\lambda)\longrightarrow(J', K', \lambda'):L\dashv
  R\longrightarrow \overline{L}\dashv\overline{R}$, $\alpha:J\longrightarrow
  J'$ and $\beta:K\longrightarrow K'$, such that they fulfill one of the
  following equivalent conditions

\begin{enumerate}

\item $\beta L\circ \lambda = \lambda' \circ \overline{L}\alpha$.

\item $\overline{R}\beta = \alpha R$. 

\end{enumerate}

\end{enumerate}

The compositions and units are the obvious ones.\\

The 2-category $\Mnd(Cat)$ is comprised of the following

\begin{enumerate}

\item The 0-cells are monads $(\mathcal{C}, S, \mu^{\Ss}, \eta^{\Ss})$. 

\item The 1-cells are pairs $(P,\varphi):(\mathcal{C}, S) \longrightarrow  (\mathcal{D}, T)$, where $P:\mathcal{C}\longrightarrow\mathcal{D}$ is a functor
  and $\varphi:TP\longrightarrow PS$ a natural transformation such that

\begin{eqnarray*}
P\mu^{\Ss}\circ \varphi S \circ T\varphi &=& \varphi\circ \mu^{\T}P\\
P\eta^{\Ss} &=& \varphi\circ \eta^{\T}P
\end{eqnarray*}

\item The 2-cells are natural transformations
  $\theta:(P,\varphi^{\Pp})\longrightarrow (Q,\varphi^{\Q})$, such that

\begin{equation*}
\varphi^{\Q}\circ T\theta = \theta S \circ \varphi^{\Pp}
\end{equation*}

\end{enumerate}

The composition and the units are the obvious ones.\\

The left 2-functor $\Phi_{\E}$ acts on $n$-cells as follows:

\begin{enumerate}

\item On 0-cells, $\Phi_{\E}(L\dashv R) = (\mathcal{C}, RL,
  R\varepsilon^{\Ll\R}L, \eta^{\R\Ll})$, \emph{i.e.} the induced monad.

\item On 1-cells, $\Phi_{\E}(J,K,\lambda) = (J,\overline{R}\lambda ): (\mathcal{C},
  RL)\longrightarrow (\mathcal{D}, \overline{R}\overline{L})$.

\item On 2-cells, $\Phi_{\E}(\alpha,\beta) = \alpha$

\end{enumerate}

$\phantom{br}$

The right 2-functor $\Psi_{\E}$ acts on $n$-cells as follows:

\begin{enumerate}

\item On 0-cells, $\Psi_{\E}(\mathcal{C}, S) = F^{\Ss}\dashv
  U^{\Ss}:\mathcal{C}\longrightarrow \mathcal{C}^{\Ss}$, \emph{i.e.} the
  Eilenberg-Moore adjunction, hence the name of the 2-adjunction after this
  action.

\item On 1-cells, $\Psi_{\E}(P, \varphi^{\Pp}) = (P, P^{\varphi}, \lambda^{\Pp}):F^{\Ss}\dashv U^{\Ss}\longrightarrow F^{\T}\dashv  U^{\T}$.

  The functor $P^{\varphi}$ acts on objects as: $P^{\varphi}(N,\chi_{\N}) = (PN,
  P\chi_{\N}\cdot \varphi^{\Pp}N)$ and the natural transformation
  $\lambda^{\Pp}$ is the mate of the commutation $U^{\T}P^{\varphi} = P
  U^{\Ss}$, that is to say

\begin{equation}\label{1510131121}
\lambda^{\Pp} = \varepsilon^{FU^{\T}}P^{\varphi}F^{\Ss}\circ F^{\T}P\eta^{UF^{\Ss}}
\end{equation}

\item On 2-cells, $\Psi_{\E}(\theta) = (\theta, \theta^{\varphi})$, where $U^{\T}\theta^{\varphi}(N,\chi_{\N}) = \theta N$.

\end{enumerate}

$\phantom{br}$

The unit of the 2-adjunction, $\eta^{\Psi\Phi_{\E}}:1_{\Adj_{\R}(Cat)}\longrightarrow \Psi_{\E}\Phi_{\E}$, has to act on 0-cells, $L\dashv R$ as follows

\begin{equation}\label{1510131134}
\eta^{\Psi\Phi_{\E}}(L\dashv R) = L\dashv R \longrightarrow F^{\R\Ll}\dashv U^{\R\Ll},
\end{equation}

$\phantom{br}$

\noindent therefore the component definition is $\eta^{\Psi\Phi_{\E}}(L\dashv R) = (1_{\mathcal{C}}, K^{\R\Ll}, \lambda^{\Ll\dashv\R})$, where
$K^{\R\Ll}:\mathcal{X}\longrightarrow \mathcal{C}^{\R\Ll}$ is \emph{the comparison functor} and $\lambda^{\Ll\dashv \R}$ is the corresponding mate, which also commutes.\\

The counit of the 2-adjunction, $\varepsilon^{\Phi\Psi_{\E}}:\Phi_{\E}\Psi_{\E}\longrightarrow 1_{\Mnd(Cat)}$, has to act on 0-cells, $(\mathcal{C}, S)$, as follows

\begin{equation}\label{1510131206}
\varepsilon^{\Phi\Psi_{\E}} (\mathcal{C}, S): (\mathcal{C}, S)\longrightarrow (\mathcal{C},S)
\end{equation}

$\phantom{br}$

\noindent therefore, the component is defined as the identity $\varepsilon^{\Phi\Psi_{\E}} (\mathcal{C}, S) = (1_{\mathcal{C}},1_{\Ss})$.

\subsection*{The Kleisli 2-adjunction for monads}

\begin{equation*}
\begin{array}{cc}
\xymatrix@C=1.6cm{
\Mnd^{\bullet}(Cat)  \ar@<-3pt>[r]_{\Psi_{\K}} & \Adj_{\Ll}(Cat)\ar@<-3pt>[l]_{\Phi_{\!\K}} 
} &\!\! ,
\end{array}
\end{equation*}

\noindent where $\Mnd^{\bullet}(Cat) \cong \Mnd(Cat^{op})$. Note the change of
direction of the 2-adjunction.\\

The 2-category $\Mnd^{\bullet}(Cat)$ is comprised of the following $n$-cells

\begin{enumerate}

\item The 0-cells are monads $(\mathcal{C}, S)$.

\item The 1-cells are pairs $(P,\psi):(\mathcal{C}, S)\longrightarrow
  (\mathcal{D}, T)$, where $P:\mathcal{C}\longrightarrow \mathcal{D}$ is a
  functor and $\psi: PS\longrightarrow TP$ is natural transformation such that
  the following requirements are fulfilled

\begin{enumerate}

\item $\mu^{\T}P\circ T\psi\circ \psi S = \psi\circ P\mu^{\Ss}$.

\item $\eta^{\T}P = \psi\circ P\eta^{\Ss}$.

\end{enumerate}

\item The 2-cells $\vartheta:(P,
  \psi^{\Pp})\longrightarrow (Q,\psi^{\Q}):(\mathcal{C}, S)\longrightarrow
  (\mathcal{D}, T)$ are natural transformations $\vartheta:P\longrightarrow
  Q:\mathcal{C}\longrightarrow \mathcal{D}$ such that 

\begin{equation*}
T\vartheta\circ \psi^{\Pp} = \psi^{\Q}\circ \vartheta S
\end{equation*}

\end{enumerate}

The 2-category $\Adj_{\Ll}(Cat)$ differs from $\Adj_{\R}(Cat)$ changing
the commutativity of the right adjoints for the left adjoints.\\

The left 2-adjunction $\Psi_{\K}$ acts on $n$-cells as follows:

\begin{enumerate}

\item On 0-cells, $\Psi_{\K}(\mathcal{C}, S) = D^{\Ss}\dashv
  V^{\Ss}:\mathcal{C}\longrightarrow\mathcal{C}_{\Ss}$, \emph{i.e.} the
  Kleisli adjunction.

\item On 1-cells, $\Psi_{\K}(P,\psi^{\Pp}) = (P, P_{\psi}, \rho^{\Pp}):D^{\Ss}\dashv  V^{\Ss}\longrightarrow D^{\T}\dashv V^{\T}$.

  The functor $P_{\psi}:\mathcal{C}_{\Ss}\longrightarrow\mathcal{D}_{\T}$ acts on
  objects as $P_{\psi}(X) = PX$ and over morphisms,
  $y^{\sharp}:X\longrightarrow Y$ in $\mathcal{C}_{\Ss}$, as $P_{\psi}(y^{\sharp}) =
  (\psi Y\cdot Py)^{\sharp}$ which makes the corresponding left adjoint square
  commuting and $\rho^{\Pp}$ is the mate for this commutative diagram.

\item On 2-cells, $\vartheta:(P,\psi^{\Pp})\longrightarrow(Q,\psi^{\Q})$,
  $\Psi_{\K}(\vartheta) = (\vartheta, \vartheta^{\psi})$, where the components
  of the natural transformation $\vartheta^{\psi}$ on objects, $X$ in
  $\mathcal{C}_{\Ss}$, is 

\begin{equation*}
\vartheta^{\psi}(X) = (\eta^{\T}QX\cdot\vartheta X)^{\sharp}
\end{equation*}

\end{enumerate}

The right adjoint $\Phi_{\K}$ only differs from the 2-functor $\Phi_{\E}$ on
how they act on the 1-cells $\Phi_{\K}(J,K,\rho) = (J, \rho L)$.\\

The unit for the 2-adjunction $\eta^{\Phi\Psi_{\K}}: 1_{\Mnd^{\bullet}(Cat)}\longrightarrow\Phi_{\K}\Psi_{\K}$  has to act on 0-cells,
$(\mathcal{C}, S)$, as follows

\begin{equation*}
\eta^{\Phi\Psi_{\K}}(\mathcal{C}, S) : (\mathcal{C}, S)\longrightarrow (\mathcal{C}, S)
\end{equation*} 

$\phantom{br}$

\noindent therefore, the component of the unit is defined as the identity
$\eta^{\Phi\Psi_{\K}}(\mathcal{C}, S) = (1_{\mathcal{C}}, 1_{\Ss})$.\\

The counit of the 2-adjunction $\varepsilon^{\Psi\Phi_{\K}}:\Psi_{\K}\Phi_{\K}\longrightarrow 1_{\Adj_{\Ll}}$ has to act on
0-cells $L\dashv R$, as follows

\begin{equation*}
\varepsilon^{\Psi\Phi_{\K}}(L\dashv R): D_{\R\Ll}\dashv V_{\R\Ll}\longrightarrow L\dashv R
\end{equation*}

$\phantom{br}$

\noindent therefore $\varepsilon^{\Psi\Phi_{\K}}(L\dashv R) = (1_{\mathcal{C}},
K_{\R\Ll}, \rho^{L\dashv R})$, where the funtor $K_{\R\Ll}$ is the \emph{Kleisli
comparison} functor. 

\subsection*{The Eilenberg-Moore 2-adjunction for Comonads}

The Eilenberg-Moore 2-adjunction for comonads 

\begin{equation*}
\begin{array}{cc}
\xymatrix@C=1.6cm{
\Adj_{\Ll}(Cat) \ar@<-3pt>[r]_{\vec{\Phi}_{\E}} & \CoMnd(Cat)\ar@<-3pt>[l]_{\vec{\Psi}_{\!\E}} 
} &\!\!\! ,
\end{array}
\end{equation*}

\noindent where 

\begin{equation*}
\Mnd(Cat^{co})\cong \CoMnd(Cat),
\end{equation*}

$\phantom{br}$

\noindent is comprised of the following structure.\\

The 2-category of comonads $\CoMnd(Cat)$ whose $n$-cells are given as follows

\begin{enumerate}

\item The 0-cells are comonads $(\mathcal{X}, G, \delta^{\G},
  \varepsilon^{\G})$.

\item The 1-cells are pairs $(P, \pi): (\mathcal{X}, G)\longrightarrow
  (\mathcal{Y}, H)$, where $P:\mathcal{X}\longrightarrow\mathcal{Y}$ is a
  functor and $\pi: PG\longrightarrow HP$ is a natural transformation such
  that 

\begin{eqnarray*}
\delta^{\Hh}P\circ \pi &=& H\pi\circ \pi G\circ P\delta^{\G}\\
\varepsilon^{\Hh}P\circ \pi &=& P\varepsilon^{\G}
\end{eqnarray*}

\item The 2-cells are natural transformations
  $\vartheta:(P,\pi^{\Pp})\longrightarrow (Q,\pi^{\Q})$, where
  $\vartheta:P\longrightarrow Q$ and it fulfills the following equation

\begin{equation*}
H\vartheta\circ\pi^{\Pp} = \pi^{\Q}\circ\vartheta G
\end{equation*}

\end{enumerate}

The compositions and units are the ones inherited from $Cat$.\\

The left 2-functor $\vec{\Phi}_{\E}$ acts on $n$-cells as follows

\begin{enumerate}

\item On 0-cells, $L\dashv R$, $\vec{\Phi}_{\E}(L\dashv R) = (\mathcal{X}, LR,
  L\eta^{\R\Ll}R, \varepsilon^{\Ll\R})$, \emph{i.e.} the induced comonad by the
  adjunction.

\item On 1-cells, $\vec{\Phi}_{\E}(J,K,\rho) = (K, \overline{L}\rho):(\mathcal{X},
  LR)\longrightarrow (\mathcal{Y}, \overline{L}\overline{R})$.

\item On 2-cells, $(\alpha, \beta):(J,K,\rho)\longrightarrow (J',K',\rho')$,
  $\vec{\Phi}_{\E}(\alpha,\beta) = \beta:K\longrightarrow K'$.

\end{enumerate} 

$\phantom{br}$

The right 2-functor $\vec{\Psi}_{\E}$ acts on $n$-cells as follows

\begin{enumerate}

\item On 0-cells, $(\mathcal{X}, G)$, $\vec{\Psi}_{\E}(\mathcal{X}, G) =
  U^{\G}\dashv F^{\G}: \mathcal{X}^{\G}\longrightarrow \mathcal{X} $, where
  $\mathcal{X}^{\G}$ is the category of Eilenberg-Moore coalgebras and the
  adjunction $U^{\G}\dashv F^{\G}$ is the usual one.

\item On 1-cells, $(P,\pi):(\mathcal{X}, G)\longrightarrow (\mathcal{Y},H)$,
  $\vec{\Psi}_{\E}(P,\pi) = (P^{\pi}, P, \rho^{\Pp}):U^{\G}\dashv
  F^{\G}\longrightarrow U^{\Hh}\dashv F^{\Hh}$. The functor
  $P^{\pi}:\mathcal{X}^{\G}\longrightarrow\mathcal{Y}^{\Hh}$ is defined over
  coalgebras $(N,\xi_{\N})$, as $P^{\pi}(N,\xi_{\N}) = (PN, \pi N\cdot
  P\xi_{\N})$. This functor makes the left square commutative and $\rho^{\Pp}$
  is the corresponding mate.

\item On 2-cells, $\vartheta:(P,\pi^{\Pp})\longrightarrow (Q,\pi^{\Q})$,
  $\vec{\Psi}_{\E}(\vartheta) = (\vartheta^{\pi},\vartheta)$, where the
  component, at $(N,\xi_{\N})$, of the natural transformation
  $\vartheta^{\pi}$ is the following

\begin{equation*}
\vartheta^{\pi}(N,\xi_{\N}) = \vartheta N
\end{equation*}

\end{enumerate}

The unit of the 2-adjunction
$\eta^{\vec{\Psi}\vec{\Phi}_{\E}}:1_{\Adj_{\Ll}(Cat)}\longrightarrow \vec{\Psi}\vec{\Phi}_{\E}$ has to act on 0-cells, $L\dashv R$, as follows

\begin{equation*}
\eta^{\vec{\Psi}\vec{\Phi}_{\E}}(L\dashv R): L\dashv R \longrightarrow
U^{\Ll\R}\dashv F^{\Ll\R}
\end{equation*}

$\phantom{br}$

\noindent therefore the component of the unit is defined as $\eta^{\vec{\Psi}\vec{\Phi}_{\E}}(L\dashv R)
= (K^{\Ll\R}, 1_{\mathcal{X}}, \rho^{\Ll\dashv\R})$, where $K^{\Ll\R}$ is the
\emph{comparison functor} for the comonad $LR$.\\

The counit of the 2-adjunction $\varepsilon^{\vec{\Phi}\vec{\Psi}_{\E}}:\vec{\Phi}_{\E}\vec{\Psi}_{\E}\longrightarrow 1_{\Mnd(Cat)}$ acts on 0-cells,
$(\mathcal{X}, G)$, as follows

\begin{equation*}
\varepsilon^{\vec{\Phi}\vec{\Psi}_{\E}}(\mathcal{X}, G): (\mathcal{X},
G)\longrightarrow (\mathcal{X}, G)
\end{equation*}

$\phantom{br}$

\noindent therefore the component of the counit is the identity in any 0-cell, that is to say,
$\varepsilon^{\vec{\Phi}\vec{\Psi}_{\E}}(\mathcal{X}, G) = (1_{\mathcal{X}}, 1_{\G})$.

\begin{center}
\section*{Acknowledgement}
\end{center}

The author would like to thank to the \emph{Consejo Nacional de Ciencia y Tecnolog\'ia}
(CONACYT) for the financial support through the grant SNI-59154.

$\phantom{br}$


\begin{thebibliography}{99}

\bibitem{agmc_gehm} AGUIAR, M. and CHASE, S.U. 
\emph{Generalized Hopf Modules for Bimonads}. Theory Appl. Categ. \textbf{27} (13), pag. 263-326. 2013.

\bibitem{eckb_setc-bejo_dila} BECK, J. Distributive Laws. In \emph{Seminar on
  Triples and Categorical Homology Theory}, B. Eckmann, Ed., vol. 80 of
  \emph{Lectures Notes in Mathematics}. Springer Berlin Heidelberg, 1969,
  pp. 119-140.

\bibitem{bral_homm} BRUGUIERES, A.; LACK, S. and VIRELIZIER, A.
\emph{Hopf Monads on Monoidal Categories}. Adv. Math. \textbf{227} (2), pag. 745-800. 2011.

\bibitem{cljs_klem} CLIMENT VIDAL, J. and SOLIVERES TUR, J.
\emph{Kleisli and Eilenberg-Moore constructions as part of a biadjoint
  situation}. Extracta Math. \textbf{25} (1), pag. 1-61. 2010.

\bibitem{lojj_reke} LOPEZ H., J.; TURCIO C., L. and VAZQUEZ-MARQUEZ, A. 
\emph{Applications of the Kleisli and Eilenberg-Moore 2-adjunctions},
2014. Available at http://arxiv.org/abs/1407.6969.

\bibitem{maep_mocg} MANES, E. and MULRY, P. 
\emph{Monad Compositions I: General constructions and recursive distributy laws}. 
Theory Appl. Categ. 18, 7
  (2007), 172-208.

\bibitem{moie_motc} MOERDIJK, I.
\emph{Monads on tensor categories}.
J. Pure Appl. Algebra \textbf{168} (2-3), pag. 189-208. 2002.

\bibitem{stro_fotm} STREET, R.
\emph{The formal theory of monads}.
J. Pure Appl. Algebra \textbf{2} (2), pag. 149-168. 1972.

\bibitem{zama_fotm} ZAWADOSKI, M.
\emph{The formal theory of monoidal monads}. 2010. Available at http://arxiv.org/abs/1012.0547

\end{thebibliography}

\end{document}